\documentclass[12pt]{article}
\usepackage{color}
\definecolor{darkblue}{rgb}{0.00,0.25,0.50}
\usepackage[colorlinks,filecolor=blue,citecolor=darkblue]{hyperref}

\usepackage{amssymb,amsmath}
\usepackage{url}
\usepackage{amscd, color}
\usepackage[english]{babel}
\begin{document}

\thispagestyle{empty}

\begin{center}
\textbf{ON THE BEST APPROXIMATION OF CERTAIN CLASSES OF PERIODIC FUNCTIONS BY
    TRIGONOMETRIC POLYNOMIALS}
\end{center}
\vskip0.5cm
\begin{center}
A.S. SERDYUK and Ie.Yu. OVSII
\end{center}
\vskip0.5cm

\abstract{ We obtain the estimates for the best approximation in the uniform metric of the
classes of $2\pi $-periodic functions whose $(\psi ,\beta )$-derivatives have a given
majorant $\omega$ of the modulus of continuity. It is shown that the estimates obtained here
are asymptotically exact under some natural conditions on the parameters $\psi ,$ $\omega $
and $\beta $ defining the classes}



\section{Introduction}
Let $L$ be the space of $2\pi $-periodic functions summable over the
period with the norm $\|f\|_1=\int_{-\pi }^{\pi }|f(t)|\,dt$
 and let $C$ be the space of $2\pi
$-periodic continuous functions $f$ with the norm
$\|f\|_C=\max\limits_t|f(t)|.$

Suppose $f\in L$ and
\begin{equation}\label{23.12.09-14:24:29}
    S[f]=\frac{a_0}{2}+\sum\limits_{
    k=1}^{\infty}(a_k\cos kx+b_k\sin kx)
\end{equation}
is its Fourier series. Suppose also that $\psi (k)$ is an arbitrary
numerical sequence and $\beta $ is a fixed real number $(\beta \in
\mathbb{R}).$ If the series
    $$
    \sum\limits_{k=1}^{\infty}\frac{1}{\psi (k)}\bigg(
    a_k\cos \Big(kx+\frac{\beta\pi}{2}\Big)+
    b_k\sin\Big(kx+\frac{\beta\pi}{2}\Big)\bigg)
    $$
is the Fourier series of a certain function $\varphi \in L,$ then $\varphi
$ is called (see, e.g., \cite{Stepanets_1987, Stepanets_2002}) the $(\psi
,\beta )$-derivative of the function $f$ and is denoted by $f^\psi _\beta
.$  The set of continuous functions $f(x)$ having $(\psi ,\beta
)$-derivatives such that $f^\psi _\beta \in H_\omega $, where
$$H_\omega =\{\varphi \in C:\
    |\varphi (t')-\varphi (t'')|\leqslant \omega (|t'-t''|)\ \ \forall
    t',t''\in \mathbb{R}\},$$
and $\omega (t)$ is a fixed modulus of continuity is usually denoted by
$C^{\psi}_{\beta}H_\omega.$

For $\psi (k)=k^{-r},\ r>0,$ the classes $C^{\psi}_{\beta}H_\omega $
 become the well-know Weyl-Nagy classes
$W^r_\beta H_\omega $  which, in turn, for $\beta =r$ coincide with the
Weyl classes $W^r_rH_\omega $ (see, e.g., \cite[Chap. 3, Sec. 4,
6]{Stepanets_2002}). For natural numbers $r$ and $\beta =r$ we obtain the
classes of periodic functions whose $r$-th derivatives are in the class
$H_\omega $.

Let $\mathfrak{M}$ be the set of all continuous functions $\psi (t)$
convex downwards for $t\geqslant 1$ and satisfying the condition
$\lim\limits_{t\to\infty}{\psi (t)}=0.$

If $\psi \in \mathfrak{M}'$, where
    $$\mathfrak{M}':=\mathfrak{M}'(\beta)=\{\psi :\
    \psi \in \mathfrak{M}\ \text{when}\ \sin\frac{\beta \pi }{2}=0\ \text{or}\ $$
    $$ \psi \in \mathfrak{M}\ \text{and}\
    \int_{1}^{\infty}\frac{\psi (t)}{t}\,dt<\infty\ \text{when}\ \sin \frac{\beta \pi }{2}\neq 0\},$$
then the classes $C^{\psi}_{\beta}H_\omega $ coincide with the classes of functions $f(x)$,
which are representable by the convolutions
    \begin{equation}\label{24.12.09-14:19:33}
    f(x)=\frac{a_0}{2}+\frac{1}{\pi }
    \int_{-\pi }^{\pi }\varphi (x+t)\Psi_\beta (t)\,dt
    ,\ \ \varphi \in H_\omega ^0,\ \ x\in \mathbb{R}
    \end{equation}
(see, e.g., \cite[p. 31]{Stepanets_1987}), where $H_\omega ^0=\{\varphi
\in H_\omega :\ \int_{-\pi }^{\pi }\varphi (t)\,dt=0\},$ and $\Psi_\beta
(t)$ is a summable function, whose Fourier series have the form
$\sum\limits_{k=1}^{\infty}\psi (k)\cos(kt+\beta \pi /2).$

The set $\mathfrak{M}$ is very inhomogeneous in the rate of convergence of
functions $\psi (t)$ to zero as $t\to\infty.$
 This is why it was suggested in
\cite[pp. 115, 116]{Stepanets_1987} (see also  \cite[Subsec.
1.3]{Stepanets and Shidlich}) to select subsets $\mathfrak{M}_0$ and
$\mathfrak{M}_C$ from $\mathfrak{M}$ as follows:
    $$\mathfrak{M}_0=\{\psi \in \mathfrak{M}:\ 0<\mu (t)\leqslant K<\infty \ \ \ \ \forall
    t\geqslant 1\},$$
    $$\mathfrak{M}_C=\{\psi \in \mathfrak{M}:\
    0<K_1\leqslant \mu (t) \leqslant K_2<\infty\ \ \
    \forall t\geqslant 1\},$$
where $\mu (t)=\mu (\psi ;t)=\frac{t}{\eta(t)-t}$,
$\eta(t)=\eta(\psi;t)=\psi^{-1}(\psi(t)/2),$
 $\psi ^{-1}(\cdot)$ is the inverse function of $\psi (\cdot)$, and
$K,$ $K_1$, $K_2$ are positive constants (possibly dependent on $\psi
 (\cdot)$). The function $\mu (\psi ;t)$ is called the modulus of half-decay
 of the function $\psi(t) .$  It is obvious that $\mathfrak{M}_C\subset\mathfrak{M}_0.$
Typical representatives of the set  $\mathfrak{M}_C$ are the functions
$t^{-r},$ $r>0,$ representatives of the set $\mathfrak{M}_0\setminus
\mathfrak{M}_C$ are the functions $\ln ^{-\alpha }(t+1),$ $\alpha >0.$ Let
$\mathfrak{M}_0'=\mathfrak{M}'\cap \mathfrak{M}_0.$ Natural
representatives of the set $\mathfrak{M}_0'$ are the functions
$\ln^{-\alpha }(t+1),$ $\alpha
>1$. It is easy to see that if $\beta =2l,$ $l\in\mathbb{Z}$, the
set $\mathfrak{M}_0'$ coincide with $\mathfrak{M}_0$. Moreover, since for
all $\psi \in \mathfrak{M}_C$
\begin{equation}\label{25.12.09-03:00:12}
    \int_{n}^{\infty}\frac{\psi (t)}{t}\,dt\leqslant K\psi (n),\ \ \ n\in
    \mathbb{N},
\end{equation}
(see \cite[p. 204]{Stepanets_2002}) then
$\mathfrak{M}_C\subset\mathfrak{M}_0'$. Throughout the paper we denote the
positive constants that may be different in different relations by $K,$
$K_i$, \mbox{$i=1,2.$}

Let us denote the best approximation of the classes
$C^{\psi}_{\beta}H_\omega$ by trigonometric polynomials $t_{n-1}(\cdot)$
of order not more than $n-1$ by $E_n(C^{\psi}_{\beta}H_\omega)$, that is
    \begin{equation}\label{13.12.09-00:56:41}
    E_n(C^{\psi}_{\beta}H_\omega)=\sup\limits_{
    f\in C^{\psi}_{\beta}H_\omega}{\inf\limits_{t_{n-1}}\|f(\cdot)-
    t_{n-1}(\cdot)\|_C}.
    \end{equation}
As is shown in \cite[p. 330]{Stepanets_1987}   if $\omega(t)$ is an
arbitrary modulus of continuity and $\psi \in \mathfrak{M}_C,$ $\beta
\in\mathbb{R}$ or $\psi \in\mathfrak{M}'_0,$ $\beta =0$, then the
following estimates hold for the quantity $E_n(C^{\psi}_{\beta}H_\omega):$
    \begin{equation}\label{9.11.09-14:47:15}
   K_1\psi (n)\omega (1/n)\leqslant
   E_n(C^{\psi}_{\beta}H_\omega)\leqslant K_2\psi (n)\omega (1/n).
    \end{equation}
When \mbox{$\psi (k)=k^{-r},$} $r>0$, $\beta\in \mathbb{R},$ the orders of
decrease of quantity (\ref{13.12.09-00:56:41}) have been known earlier
\cite{EFIMOV61} (see also  \cite[p. 508]{TIMAN-60}).

It should be noted that unlike order estimates, exact values for the
quantity $E_n(C^{\psi}_{\beta}H_\omega)$ have been found for \mbox{$\psi
(k)=k^{-r},$} \mbox{$r\in \mathbb{N}$,} $\beta =r$ and for the convex
upwards modulus of continuity by Korneichuk \cite{KORNEICHUK-3.5} (see
also \cite[p. 319]{KORNEICHUK-5}, \cite[p. 344]{DeVore_Lorentz}). The
similar problem on the class of real-valued functions defined on the
entire real axis and having the $r$-th continuous derivatives $f^{(r)}$
such that $\omega (f^{(r)};t)\leqslant \omega (t),$ $t\in [0,\infty)$, is
solved in the paper of Ganzburg \cite{Ganzburg}.

The aim of the present work is to study the rate of decrease of quantity
(\ref{13.12.09-00:56:41})  when $\psi \in\mathfrak{M}'_0$ and $\beta\in
\mathbb{R}.$


\section{Main Results}
The following statements are true.

\hypertarget{14.04.11-19:15:10}{\textbf{Theorem {1}.}} \emph{Let $\psi \in\mathfrak{M}'_0$,
$\beta\in \mathbb{R}$ and let $\omega (t)$ be an arbitrary modulus of continuity. Then, as
$n\to\infty,$
\begin{equation}\label{19.12.09-20:55:52}
    E_n(C^{\psi}_{\beta}H_\omega)=\frac{\theta _n(\omega )}{\pi }
    \big|\sin \frac{\beta \pi }{2}\big|\int_{
    0}^{1/n}\psi \Big(\frac{1}{t}\Big)\frac{\omega (t)}{t}\,dt+O(1)
    \psi (n)\omega (1/n),
\end{equation}
where $\theta _n(\omega )\in[2/3,1]$ and $O(1)$ is a quantity uniformly bounded in $n$ and
$\beta .$ If $\omega (t)$ is a convex upwards modulus of continuity, then $\theta _n(\omega
)=1.$}

We give an example of functions $\psi $ and $\omega $ for which
(\ref{19.12.09-20:55:52}) is an asymptotic formula.

\hypertarget{14.04.11-19:32:42}{\textbf{{Example 1.}}} Let \mbox{$\psi (t)=\ln^{-\gamma
}(t+1),$ $\gamma
>1$, $\beta
\neq 2l,$ $l\in\mathbb{Z}$ }  and
    $$\omega (t)=
   \begin{cases}
    0, &  \ \ \ t=0, \\
    \ln^{-\alpha }\left(\frac{1}{t}+1\right)
    , & \ \ \
    t>0,\ \ \ 0<\alpha \leqslant 1.
  \end{cases}$$
Then by virtue of (\ref{19.12.09-20:55:52}) the following asymptotic
formula holds as $n\to\infty$:
    $$E_n(C^{\psi}_{\beta}H_\omega)=
    \ln^{-(\gamma +\alpha) }(n+1)\bigg(\frac{1}{
    \pi (\gamma +\alpha-1)}\big|\sin\frac{\beta \pi }{2}\big|\ln n+O(1)\bigg),$$
where $O(1)$ is a quantity uniformly bounded in $n$ and $\beta .$
\\

Note that if
\begin{equation}\label{17.08.10-17:16:09}
    \lim\limits_{n\to\infty}{\frac{|\psi '(n)|n}{\psi (n)}}=0,\ \
    \psi '(n):=\psi '(n+),
\end{equation}
and
\begin{equation}\label{17.08.10-17:13:36}
    \lim\limits_{n\to\infty}{\frac{\omega' (1/n)}{
    \omega (1/n)n}}=0,\ \ \omega '(1/n):=\omega '(1/n+),
\end{equation}
then equalities
    $$\lim\limits_{n\to\infty}{\frac{
    \psi (n)\omega (1/n)}{\int_{
    0}^{1/n}\psi (\frac{1}{t})\frac{\omega (t)}{t}\,dt}}=
    \lim\limits_{n\to\infty}\frac{|\psi '(n)|n}{\psi (n)}+
    \lim\limits_{n\to\infty}{\frac{\omega' (1/n)}{
    \omega (1/n)n}}=0,$$
are valid.

Therefore from Theorem \hyperlink{14.04.11-19:15:10}{1} we obtain.

\textbf{Corollary 1.} \emph{Assume that $\psi \in \mathfrak{M}_0'$, $\beta \neq 2l,$ $l\in
\mathbb{Z}$, $\omega (t)$ is a convex upwards modulus of continuity and conditions
$(\ref{17.08.10-17:16:09})$ and $(\ref{17.08.10-17:13:36})$ are fulfilled. Then the
following asymptotic formula holds as $n\to\infty:$
    $$E_n(C^{\psi}_{\beta}H_\omega)=\frac{1}{\pi}
    \big|\sin \frac{\beta \pi }{2}\big|\int_{
    0}^{1/n}\psi \Big(\frac{1}{t}\Big)\frac{\omega (t)}{t}\,dt+O(1)
    \psi (n)\omega (1/n),$$
where $O(1)$ is a quantity uniformly bounded in $n$ and $\beta .$}

The functions $\psi $ and $\omega $ from Example \hyperlink{14.04.11-19:32:42}{1} can serve
as an example of the functions which satisfy conditions $(\ref{17.08.10-17:16:09})$ and
$(\ref{17.08.10-17:13:36})$, respectively.

Relation (\ref{19.12.09-20:55:52}) implies that if $\psi \in
\mathfrak{M}_0'$ and
\begin{equation}\label{1.12.09-00:38:51}
        \big|\sin\frac{\beta \pi }{2}\big|
        \int_{0}^{1/n}\frac{\omega (t)}{t}\,dt=O(1)\omega (1/n),\ \ \beta
        \in \mathbb{R},
   \end{equation}
or $\psi \in\mathfrak{M}_C$ (see (\ref{25.12.09-03:00:12})), then 
    $$E_n(C^{\psi}_{\beta}H_\omega)=O(1)\psi (n)\omega (1/n).$$
Taking into account that function $\psi (t)$ is monotonically decreasing
for $t\geqslant 1$ and using the estimate
\begin{equation}\label{24.12.09-19:54:08}
    E_n(C^{\psi}_{\beta}H_\omega)\geqslant K\psi (n)\omega (1/n)\ \ \ \forall
    \psi \in \mathfrak{M}',
\end{equation}
(see \cite[pp. 329, 330]{Stepanets_1987}), by virtue of relation
(\ref{19.12.09-20:55:52}) we arrive at the following statement:

\textbf{Corollary 2.} \emph{Let $\beta \in \mathbb{R}$ and let $\omega (t)$ be an arbitrary
modulus of continuity. If $\psi \in \mathfrak{M}_C$ or $\psi \in \mathfrak{M}_0'$ and
$\omega (t)$ satisfies condition $(\ref{1.12.09-00:38:51})$, then
\begin{equation}\label{1.12.09-00:51:42}
K_1\psi (n)\omega (1/n)\leqslant
   E_n(C^{\psi}_{\beta}H_\omega)\leqslant K_2\psi (n)\omega (1/n),
\end{equation}
where $K_1$ and $K_2$ are positive constants.}

Thus, estimates (\ref{9.11.09-14:47:15}) obtained by Stepanets \cite[p.
330]{Stepanets_1987} (see also \cite[Chap. 5, Sec. 22; Chap. 7, Sec.
4]{Stepanets_2002}) for the arbitrary modulus of continuity $\omega (t)$
and for $\psi \in \mathfrak{M}_C,$ $\beta \in\mathbb{R}$ or for $\psi
\in\mathfrak{M}'_0,$ $\beta =0,$ hold also in the case when $\psi \in
\mathfrak{M}_0'$, $\beta \neq 0$ and $\omega (t)$ satisfies condition
(\ref{1.12.09-00:38:51}). For example, the function $\omega (t)=t^\alpha
,$ \mbox{$0<\alpha\leqslant 1,$} satisfies (\ref{1.12.09-00:38:51}).

\section{Proof of Theorem 1}

Suppose that all conditions of the theorem are satisfied. Let us carry out
the proof in two stages.

\textbf{1.} We shall find an upper estimate for
 $E_n(C^{\psi}_{\beta}H_\omega).$

We set
    \begin{equation}\label{23.12.09-14:37:17}
    U_{n-1}^\psi (f;x)=\frac{a_0}{2}+\sum\limits_{
    k=1}^{n-1}\Big(1-\frac{\psi (n)}{\psi (k)}\frac{k^2}{n^2}\Big)(
    a_k\cos kx+b_k\sin kx),\ \ n\in\mathbb{N},
    \end{equation}
where $a_k$ and $b_k$ are the Fourier coefficients of a function $f\in
C^{\psi}_{\beta}H_\omega.$ Show that for the quantity
    $$\mathcal{E}_n(C^{\psi}_{\beta}H_\omega)=\sup\limits_{
    f\in C^{\psi}_{\beta}H_\omega}{\|f(\cdot)-U_{n-1}^\psi (f;\cdot)\|_C}$$
the inequality
    \begin{equation}\label{23.12.09-13:55:31}
    \mathcal{E}_n(C^{\psi}_{\beta}H_\omega)\leqslant
    \frac{1}{\pi }\big|\sin\frac{\beta \pi }{2}\big|\int_{0}^{1/n}\psi \Big(
    \frac{1}{t}\Big)\frac{\omega (t)}{t}\,dt+O(1)\psi (n)\omega (1/n),
    \end{equation}
is true. Since
\begin{equation}\label{23.12.09-13:55:02}
    E_n(C^{\psi}_{\beta}H_\omega)\leqslant\mathcal{E}_n(
    C^{\psi}_{\beta}H_\omega),
\end{equation}
then the required upper estimate for $E_n(C^{\psi}_{\beta}H_\omega)$
follows from (\ref{23.12.09-13:55:31}).

For further reasoning, we need the one statement, which follows from the
book \cite[p. 65]{Stepanets_1987}. We will give a few notations before
formulating it. Let $f$ be a summable function, whose Fourier series have
the form (\ref{23.12.09-14:24:29}). Further, let $\lambda _n=\{\lambda
_1(u),\lambda _2(u),\ldots,\lambda _n(u)\}$ be a collection of continuous
functions on $[0,1]$ such that $\lambda (k/n)= \lambda _k^{(n)},$
$k=\overline{0,n},$ $n\in\mathbb{N},$ where $\lambda _k^{(n)}$ are
elements of the triangular matrix $\Lambda=\|\lambda _k^{(n)}\|$,
$k=\overline{1,n},$ $\lambda _0^{(n)}=1,$ that determine a polynomial of
the form
    \begin{equation}\label{24.12.09-14:09:14}
    U_n(f;x;\Lambda )=
    \frac{a_0}{2}+\sum\limits_{
    k=1}^{n}\lambda _k^{(n)}(a_k\cos kx+b_k\sin kx),\ \ n\in \mathbb{N}.
    \end{equation}

The following statement is true:

\hypertarget{14.04.11-19:36:38}{\textbf{Lemma A}} \cite[p. 65]{Stepanets_1987}.
\emph{Suppose that $f\in C^{\psi}_{\beta}H_\omega $ and $\tau _n(u)$ is the continuous
function defined by relation}
     \emph{\begin{equation}\label{2/06-5}
        \tau  _n(u)=\tau  _n(u;\lambda ;\psi )=
 \begin{cases}
    (1-\lambda _n(u))\psi (1)nu, & \text{\ \ \ \
    $0\leqslant u\leqslant \frac{1}{n},$}\\
    (1-\lambda _n(u))\psi (nu), & \text{\ \ \ \
    $\frac{1}{n}\leqslant u\leqslant 1,$}\\
    \psi (nu), & \text{\ \ \ \ $u\geqslant 1,$}
  \end{cases}
  \end{equation}
and such that its Fourier transform
    $$\widehat{\tau }_n(t):=\widehat{\tau }_n(t;\beta )=\frac{1}{\pi }\int_{0}^{\infty}
    \tau_n (u)\cos\Big(ut+\frac{\beta  \pi }{2}\Big)
    \,du,\ \ \ \beta   \in\mathbb{R},$$
is summable on the whole real line, i.e.
    $\int_{-\infty}^{\infty}|\widehat{\tau }_n(t)|\,dt<\infty.$
Then at any point x the following equality holds}:
    \begin{equation}\label{2/06-7}
        f(x)-U_{n}(f;x;\Lambda  )=\int_{-
        \infty}^{\infty}f^\psi _\beta \Big(x+\frac{t}{n}
        \Big)\widehat{\tau }_n(t)\,dt,\ \ \ n\in\mathbb{N}.
    \end{equation}

Using Lemma \hyperlink{14.04.11-19:36:38}{A}, let us show that
\begin{equation}\label{23.12.09-14:31:56}
    f(x)-U_{n-1}^\psi (f;x)=\int_{
    -\infty}^{\infty}f^\psi _\beta \Big(x+\frac{t}{n}\Big)\widehat{\tau }
    _n(t)\,dt\ \ \ \forall f\in
C^{\psi}_{\beta}H_\omega,\ \ n\in \mathbb{N},
\end{equation}
where $\widehat{\tau }_n(t)$ is the Fourier transform of the function
\begin{equation}\label{23.12.09-14:35:01}
  \tau _n(u)=\tau _n(u;\psi )=\begin{cases}
    \psi (n)u^2, & 0\leqslant u\leqslant 1, \\
    \psi (nu), & u\geqslant 1.
  \end{cases}
\end{equation}

Since polynomial (\ref{23.12.09-14:37:17}) can be represented in the form
    $$U_{n-1}^\psi (f;x)=
    \frac{a_0}{2}+\sum\limits_{
    k=1}^{n}\lambda ^\psi (k/n)(a_k\cos kx+b_k\sin kx),$$
where $\lambda ^\psi (k/n)$ are the values of continuous function
    \begin{equation}\label{24.12.09-14:52:56}
    \lambda ^\psi (u)=\lambda ^\psi_n (u)=
  \begin{cases}
    1-\frac{\psi (n)}{\psi (1)}\frac{u}{n}, & 0\leqslant u\leqslant \frac{1}{n}, \\
    1-\frac{\psi (n)}{\psi (nu)}u^2, & \frac{1}{n}\leqslant u\leqslant 1
  \end{cases}
    \end{equation}
at the points $u=k/n$ and
    $$\tau _n(u)=\tau _n(u;\psi )=
  \begin{cases}
    (1-\lambda^\psi (u))\psi (1)nu, & 0\leqslant
    u\leqslant \frac{1}{n}, \\
    (1-\lambda ^\psi (u))\psi (nu), & \frac{1}{n}\leqslant u\leqslant 1, \\
    \psi (nu), & u\geqslant 1,
  \end{cases}
    $$
then it follows from Lemma \hyperlink{14.04.11-19:36:38}{A} that for proving
(\ref{23.12.09-14:31:56}) it is sufficient to establish the inequality
\begin{equation}\label{23.12.09-15:02:09}
    \int_{-\infty}^{\infty}|\widehat{\tau }_n(t)|\,dt<\infty.
\end{equation}
With this aim we put
    $$\mu _n(u)=
  \begin{cases}
    \psi (n)(u^2-u), & 0\leqslant u\leqslant 1, \\
    0, & u\geqslant 1,
  \end{cases}\ \ \ \ \ \nu _n(u)=\tau _n(u)-\mu _n(u).
    $$
Integrating twice by parts, we get
    $$\widehat{\mu }_n(t):=\widehat{\mu }_n(t;\beta )=\frac{1}{\pi }\int_{0}^{\infty}\mu _n(u)\cos\Big(ut+\frac{\beta \pi }{2}\Big)\,du=\frac{
    O(1)}{t^2},\ \ t>0,$$
which yields
\begin{equation}\label{23.12.09-15:12:06}
    \int_{|t|\geqslant 1}|\widehat{\mu }_n(t)|\,dt<\infty.
\end{equation}
It is obvious that
\begin{equation}\label{23.12.09-15:14:31}
    \int_{|t|\leqslant 1}|\widehat{\mu }_n(t)|\,dt<\infty.
\end{equation}
Taking (\ref{23.12.09-15:12:06}), (\ref{23.12.09-15:14:31}) together and
using the estimates
    $$\int_{-\infty}^{\infty}|\widehat{\nu  }_n(t)|\,dt<\infty
    \ \ \ \forall \psi \in \mathfrak{M}_0'$$
(see, e.g., \cite[p. 174]{Stepanets_2002}) and
    $$|\widehat{\tau }_n(t)|\leqslant |\widehat{\mu }_n(t)|+|\widehat{
    \nu }_n(t)|,$$
we obtain (\ref{23.12.09-15:02:09}).

Furthermore, since the function $\tau _n(u)$ satisfies all conditions of
Lemma 3 from  \cite{TELYKOVSKIY-61} according to which
    $$\tau _n(u)=\int_{-\infty}^{\infty}
    \cos \Big(ut+\frac{\beta \pi }{2}\Big)\widehat{\tau }_n(t)\,dt,\ \ \ u\geqslant 0,$$
we have
    $$\int_{-\infty}^{\infty}\widehat{\tau }_n(t)\,dt=\frac{\tau _n(0)}{\cos
    \frac{\beta \pi }{2}}=0,\ \ \ \beta \neq 2l-1,\ \ l\in \mathbb{Z}.$$
If $\beta =2l-1,$ $l\in\mathbb{Z},$ the equality $
\int_{-\infty}^{\infty}\widehat{\tau }_n(t)\,dt=0$ is obvious, because
$\widehat{\tau }_n(t)$ is odd. Hence, starting from
(\ref{23.12.09-14:31:56}) we can write
\begin{equation}\label{23.12.09-15:37:11}
    f(x)-U_{n-1}^\psi (f;x)=\int_{
    -\infty}^{\infty}\bigg(f^\psi _\beta \Big(x+\frac{t}{n}\Big)-
    f^\psi _\beta (x)\bigg)\widehat{\tau }
    _n(t)\,dt\ \ \ \forall f\in
C^{\psi}_{\beta}H_\omega,\ \ n\in \mathbb{N}.
\end{equation}
Since $f^\psi _\beta \in H_\omega ^0$ and, as it is not hard to see, for
every $\varphi \in H_\omega^0$ function $\varphi _1(u)=\varphi (u+h),$
$h\in \mathbb{R},$ also belongs to $H_\omega ^0,$ then using the notation
    $$\delta (t,\varphi )=\varphi (t)-\varphi (0),$$
it follows from (\ref{23.12.09-15:37:11}) that
\begin{equation}\label{23.12.09-15:42:55}
    \mathcal{E}_n(C^{\psi}_{\beta}H_\omega)\leqslant
    \sup\limits_{\varphi \in H_\omega ^0}{\bigg|\int_{
    -\infty}^{\infty}\Big(\varphi \Big(\frac{t}{n}\Big)-\varphi (0
    )\Big)\widehat{\tau }_n(t)\,dt\bigg|}$$
    $$=\sup\limits_{\varphi \in H_\omega ^0}{\bigg|\int_{
    -\infty}^{\infty}\delta \Big(\frac{t}{n},\varphi \Big)\widehat{\tau
    }_n(t)\,dt\bigg|.}
\end{equation}
Now we shall simplify the integral in the right-hand side of
(\ref{23.12.09-15:42:55}) without loss of its principal value. The
following relations are true:
\begin{equation}\label{15/07/09-15:19:31}
    \int_{-\infty}^{\infty}\delta\Big(\frac{t}{n},\varphi\Big)
    \widehat{\tau }_n(t)\,dt$$
    $$=\frac{\cos\frac{\beta \pi}{2}}{\pi }\int_{-\infty}^{
    \infty}\delta\Big(\frac{t}{n},\varphi\Big) \int_{0}^{\infty}\tau _n(u)
    \cos ut\,du\,dt$$
    $$-\frac{\sin\frac{\beta \pi}{2}}{\pi }\int_{-\infty}^{
    \infty}\delta\Big(\frac{t}{n},\varphi\Big) \int_{0}^{\infty}\tau _n(u)
    \sin ut\,du\,dt$$
    $$=\frac{\cos\frac{\beta \pi}{2}}{\pi }\int_{-\infty}^{
    \infty}\delta\Big(\frac{t}{n},\varphi\Big) \int_{0}^{\infty}\tau _n(u)
    \cos ut\,du\,dt$$
    $$-\frac{\sin \frac{\beta \pi}{2}}{\pi }\Bigg(\int_{|t|\geqslant 1
    }\delta\Big(\frac{t}{n},\varphi\Big) \int_{
    0}^{\infty}\tau _n(u)\sin ut\,du\,dt$$
    $$+\int_{|t|\leqslant 1}\delta\Big(\frac{t}{n},\varphi\Big) \int_{
    0}^{1}\tau _n(u)\sin ut\,du\,dt$$
    $$+\int_{|
    t|\leqslant 1}\delta\Big(\frac{t}{n},\varphi\Big)
    \int_{1}^{\infty}\psi (nu)\sin ut\,du\,dt\Bigg).
\end{equation}
Integrating by parts, taking into account the equality $\tau _n(0)=\tau
_n(\infty)=0$ and assuming that \mbox{$\psi '(u):=\psi '(u+),$} we have
\begin{equation}\label{15/07/09-15:30:03}
    \int_{0}^{\infty}\tau _n(u)\cos ut\,du=-\frac{1}{t}\int_{
    0}^{\infty}\tau _n'(u)\sin ut\,du$$
    $$=-\frac{2\psi (n)}{t}\int_{0}^{1}u\sin ut\,du-\frac{n}{t}\int_{1}^{\infty}
    \psi '(nu)\sin ut\,du
\end{equation}
and similarly,
    \begin{equation}\label{15/07/09-15:34:51}
        \int_{0}^{\infty}\tau _n(u)\sin ut\,du
        =\frac{2\psi (n)}{t}\int_{0}^{1}u\cos ut\,du+\frac{n}{t}\int_{1}^{\infty
        }\psi '(nu)\cos ut\,du.
    \end{equation}
Combining (\ref{15/07/09-15:19:31})--(\ref{15/07/09-15:34:51}), we obtain
\begin{equation}\label{23.12.09-17:11:58}
    \int_{-\infty}^{\infty}\delta\Big(\frac{t}{n},\varphi\Big)
    \widehat{\tau }_n(t)\,dt$$
    $$=-\frac{\sin\frac{\beta \pi }{2}}{\pi }\int_{|
    t|\leqslant 1}\delta\Big(\frac{t}{n},\varphi\Big)
    \int_{1}^{\infty}\psi (nu)\sin ut\,du\,dt+r_n(\psi ,\varphi ,\beta ),
    \ \varphi \in H_\omega ^0,\ n\in \mathbb{N},
\end{equation}
where
\begin{equation}\label{23.12.09-17:15:29}
    r_n(\psi ,\varphi ,\beta )=$$
    $$=\frac{
    \cos\frac{\beta \pi }{2}}{\pi }\Bigg(
    -2\psi (n)\int_{-\infty}^{\infty}\delta\Big(\frac{t}{n},\varphi\Big)
    \frac{1}{t}\int_{0}^{1}u\sin ut\,du\,dt$$
    $$-n\int_{-\infty}^{\infty}\delta\Big(\frac{t}{n},\varphi\Big)
    \frac{1}{t}\int_{1}^{\infty}\psi '(nu)\sin ut\,du\,dt\Bigg)$$
    $$-\frac{\sin\frac{\beta \pi }{2}}{\pi }\Bigg(
    2\psi (n)\int_{|t|\geqslant 1}\delta\Big(\frac{t}{n},\varphi\Big) \frac{1}{t}\int_{
    0}^{1}u\cos ut\,du\,dt$$
    $$+n\int_{|t|\geqslant 1}\delta\Big(\frac{t}{n},\varphi\Big)
    \frac{1}{t}\int_{1}^{\infty}\psi '(nu)\cos ut\,du\,dt$$
    $$+\int_{|t|\leqslant 1}\delta\Big(\frac{t}{n},\varphi\Big)
    \int_{0}^{1}\tau _n(u)\sin ut\,du\,dt\Bigg)
    =:\frac{\cos\frac{\beta \pi }{2}}{\pi }\sum\limits_{i=1}^{2}J_{i,n}-
    \frac{\sin\frac{\beta \pi }{2}}{\pi }\sum\limits_{
    i=3}^{5}J_{i,n}.
\end{equation}
Let us show that
\begin{equation}\label{23.12.09-17:29:28}
    r_n(\psi ,\varphi ,\beta )=O(1)\psi (n)\omega (1/n).
\end{equation}
Since for $t\in[-1,1]$ the quantity
    $$\frac{1}{t}\int_{0}^{1}u\sin ut\,du$$
is bounded by a constant, then using the inequality $|\delta (t,\varphi
)|\leqslant \omega (|t|),$ we get
\begin{equation}\label{23.12.09-17:34:42}
    J_{1,n}=-2\psi (n)\int_{
    |t|\geqslant 1}\delta\Big(\frac{t}{n},\varphi\Big) \frac{1}{t}
    \int_{0}^{1}u\sin ut\,du\,dt+O(1)\psi (n)\omega (1/n).
\end{equation}
To estimate the integral in (\ref{23.12.09-17:34:42}) we establish the
following auxiliary statements.

\hypertarget{18.08.10-16:16:42}{\textbf{Lemma 1.}} \emph{On every interval
$(x_k^{(i)},x_{k+1}^{(i)}),\ x_k^{(i)}=(2k-1+i)\pi /2a,$ \mbox{$i=0,1,$} \mbox{$k\in
\mathbb{N}$,} $a>0,$ the function
    $$\int_{x}^{\infty}
    \frac{1}{t}\int_{0}^{a}u^s\sin \Big(ut+\frac{i\pi }{2}\Big)\,du\,dt,\ \ \
    x>0,\ \ \ s\geqslant 1,$$
has at least one zero.}

\emph{Proof.} We will give a proof of the lemma only for the case $i=0$, because the proof
in case $i=1$ is similar. On the basis of the estimate
    $\big|\int_{x}^{\infty}
    \frac{\sin t}{t}\,dt\big|\leqslant \frac{2}{x},$ $x>0$
(see, e.g., \cite[p. 5]{Tabl}, \cite[p. 191]{Stepanets-77}) it is simple
to see that the integral
    $$\int_{x}^{\infty}\frac{u^s\sin ut}{t}\,dt=
    u^s\int_{ux}^{\infty}\frac{\sin t}{t}\,dt$$
converges uniformly with respect to $u\in[0,a],\ a>0$. Therefore, changing
the order of integration, we obtain
    $$S(x):=\int_{x}^{\infty}
    \frac{1}{t}\int_{0}^{a}u^s\sin ut\,du\,dt=\int_{0}^{a}u^s
    \int_{x}^{\infty}\frac{\sin ut}{t}\,dt\,du.$$
Making the change of variables and integrating by parts, we have
    $$S(x)=\int_{0}^{a}
    u^s\int_{ux}^{\infty}\frac{\sin t}{
    t}\,dt\,du$$
    $$=\frac{1}{s+1}\Bigg(a^{s+1}
    \int_{ax}^{\infty}\frac{
    \sin t}{t}\,dt+\int_{0}^{a}u^s\sin ux\,du\Bigg)$$
    $$=\frac{1}{s+1}
    \Bigg(a^{s+1}\int_{ax}^{\infty}\frac{\sin t}{t}\,dt
    -a^s\frac{\cos ax}{x}+\frac{s}{x}\int_{0}^{a}
    u^{s-1}\cos ux\,du\Bigg).$$
Hence, taking into account the equation
    $$\int_{ax}^{\infty}
    \frac{\sin t}{t}\,dt=\frac{\cos ax}{ax}-\int_{ax
    }^{\infty}\frac{\cos t}{t^2}\,dt$$
we get
    \begin{equation}\label{1/04-94}
        S(x)=
        \frac{1}{s+1}\Bigg(
        -a^{s+1}\int_{ax}^{\infty}\frac{\cos t}{t^2}\,dt+
        \frac{s}{x^{s+1}}\int_{0}^{ax}
        u^{s-1}\cos u\,du\Bigg).
    \end{equation}
On every interval $(t_j,t_{j+1}),\ t_j=(2j+1)\pi /2,$
\mbox{$j=0,1,\ldots,$} the function
    $\int_{x}^{\infty}\frac{\cos t}{t^2}\,dt$
vanishes with a change of sign at some point $\widetilde{x}_j$. Since
    $$\int_{\pi /2}^{\infty}\frac{\cos t}{t^2
    }\,dt=-\int_{\pi /2}^{\infty}
    \frac{\sin t}{t}\,dt<0,$$
then for any $k\in \mathbb{N}$
\begin{equation}\label{1/04-95}
    \text{sign}\int_{(2k-1)\pi /2}^{\infty
    }\frac{\cos t}{t^2}\,dt=(-1)^{k}.
\end{equation}
Further,  we have
    $$\int_{0}^{(2k-1)\pi /2}
    u^{s-1}\cos u\,du=\alpha _0+\sum\limits_{
    j=1}^{k-1}\alpha _j,$$
where
    $$\alpha _0=\int_{0}^{\pi /2}
    u^{s-1}\cos u\,du,\ \ \ \ \ \alpha _j=\int_{(2j-1)\pi /2}^{
    (2j+1)\pi /2}
    u^{s-1}\cos u\,du.$$
If $k=1$, then
    \begin{equation}\label{14.01.10-19:14:55}
    \text{sign}\int_{0}^{(2k-1)\pi /2}
    u^{s-1}\cos u\,du=\text{sign}\,\alpha _0=1.
    \end{equation}
Let $k=2,3,\dots$ Since the function $u^{s-1}$ does not decrease
$(s\geqslant 1)$ for $u\geqslant 0$, we can write
      $$|\alpha _0|<|\alpha _{j}|\leqslant |\alpha _{j+1}|,\ \ \
    \ \ \ \ j\geqslant 1,$$
and respectively
    \begin{equation}\label{7/08-1}
    \text{sign}\int_{0}^{(2k-1)\pi /2}
    u^{s-1}\cos u\,du$$
    $$=
    \text{sign}\int_{(2k-3)\pi /2}^{
    (2k-1)\pi /2}
    u^{s-1}\cos u\,du=(-1)^{k+1},\ \ k=2,3,\ldots
    \end{equation}
Taking account of (\ref{1/04-94})--(\ref{7/08-1}), we have
    \begin{equation}\label{11.02.09-03:42:24}
    \text{sign}\,S\bigg(\frac{
    2k-1}{2a}\pi \bigg)=(-1)^{k+1},\ \ \ k\in \mathbb{N},\ \ \
    a>0.
    \end{equation}
The function $S(x)$ is continuous for any $x>0.$ Therefore, it follows from
(\ref{11.02.09-03:42:24}) that on every interval $(x_k,x_{k+1})$, where \mbox{$x_k=(2k-1)\pi
/2a,\ k\in\mathbb{N},\ a>0,$} the function $S(x)$ has the required zero. Lemma
\hyperlink{18.08.10-16:16:42}{1} is proved.\hspace{\stretch{1}}$\blacksquare$

\hypertarget{18.08.10-16:25:32}{\textbf{Lemma 2.}} \emph{Let $\varphi  \in H_\omega $,
$1\leqslant a\leqslant n,\ n\in\mathbb{N}$ and $s\geqslant 1$. Then for $i=0,1,$ the
following estimate holds$:$
\begin{equation}\label{7/08-3}
    \int_{
    |t|\geqslant 1}\bigg(\varphi \Big(\frac{t}{n}\Big)-\varphi (0)\bigg)
    \frac{1}{t}\int_{
    0}^{a/n}u^s\sin\Big(ut+\frac{i\pi }{2}\Big)\,du\,dt=O(1
    )\omega (1/n),
\end{equation}
where
 $O(1)$ is a quantity uniformly bounded in $n,$ $\varphi ,$ $a$ and $s.$}

\emph{Proof.} Making the change of variables, we get
\begin{equation}\label{11/08-10}
    \int_{
    |t|\geqslant 1}\bigg(\varphi \Big(\frac{t}{n}\Big)-\varphi
    (0)\bigg)\frac{1}{t}\int_{
    0}^{a/n}u^s\sin\Big(ut+\frac{i\pi }{2}\Big)\,du\,dt$$
    $$=\frac{1}{n^{s+1}}\int_{|t|\geqslant 1/n}
    (\varphi (t)-\varphi (0))\frac{1}{t}\int_{0}^{a}u^s\sin\Big(ut+
    \frac{i\pi }{2}\Big)\,du\,dt, \ \ \ i=0,1.
\end{equation}
Let us denote by $t_k^{(i)}$ the zero of function
    $$\int_{x}^{\infty}\frac{1}{t}
    \int_{0}^{a}u^s\sin\Big(ut+\frac{i\pi
    }{2}\Big)\,du\,dt, \ \ \ i=0,1,$$
on interval $(x_k^{(i)},x_{k+1}^{(i)})$, $x_k^{(i)}=\frac{2k-1+i}{2a}\pi,$ which exists
according to Lemma \hyperlink{18.08.10-16:16:42}{1}. Using the notation $\delta (t)=\varphi
(t)-\varphi (0)$, we have
\begin{equation}\label{7/08-4}
    \bigg|\int_{1/n}^{\infty}\delta (t)\frac{1}{t}
    \int_{0}^{a}u^s\sin\Big(ut+\frac{i\pi }{2}\Big)\,du\,dt\bigg|
    $$
    $$=  \bigg|{\int_{1/n}^{{t}_1^{(i)}}
    \delta (t)\frac{1}{t}
    \int_{0}^{a}u^s\sin\Big(ut+\frac{i\pi }{2}\Big)
    \,du\,dt}$$
    $$+\sum\limits_{k =1}^{\infty}\int_{{t}_k ^{(i)}}^{{t}_{k+1} ^{(i)}}
    \big(\delta (t)-\delta (t_k^{(i)})\big)\frac{1}{t}\int_{0}^{a}u^s\sin\Big(ut+\frac{i\pi }{2}\Big)
    \,du\,dt\bigg|$$
    $$\leqslant  {
    \omega  ({t}_1^{(i)})\int_{1/n}^{{t}_1^{(i)}}\frac{1}{t}
    \bigg|\int_{0}^{a}u^s\sin\Big(ut+\frac{i\pi }{2}\Big)
    \,du\bigg|\,dt}$$
    $$+\omega (\Delta_i)\int_{{t}_1^{(i)}}^{\infty}
    \frac{1}{t}\bigg|\int_{0}^{a}u^s\sin\Big(ut+\frac{i\pi }{2}\Big)
    \,du\bigg|\,dt,
\end{equation}
where $\Delta_i =\sup\limits_{k}{(t_{k+1}^{(i)}-t_k^{(i)})}.$ Since
${t}_1^{(i)}<\frac{2\pi }{a} $ and $\Delta_i <\frac{2\pi }{a},$ it follows
from (\ref{7/08-4}) that
\begin{equation}\label{7/08-5}
    \bigg|\int_{1/n}^{\infty}\delta (t)\frac{1}{t}
    \int_{0}^{a}u^s\sin\Big(ut+\frac{i\pi }{2}\Big)\,du\,dt\bigg|
      $$
    $$<\omega \Big(\frac{2\pi }{a}\Big)\int_{1/n}^{\infty}
    \frac{1}{t}\bigg|\int_{0}^{a}u^s\sin\Big(ut+\frac{i\pi }{2}\Big)
    \,du\bigg|\,dt.
\end{equation}
After integrating by parts it is easy to see, that
    \begin{equation}\label{11.02.09-03:55:57}
    \bigg|\int_{0}^{a}u^s\sin\Big(ut+\frac{i\pi }{2}\Big)\,du\bigg|\leqslant
    \frac{2a^s}{t},\ \ t>0,\ \ i=0,1.
    \end{equation}
From (\ref{7/08-5}) and (\ref{11.02.09-03:55:57}) follows the inequality
\begin{equation}\label{7/08-6}
    \bigg|\int_{1/n}^{\infty}\delta (t)\frac{1}{t}
    \int_{0}^{a}u^s\sin\Big(ut+\frac{i\pi }{2}\Big)\,du\,dt
    \bigg|
    < 2a^s\omega \Big( \frac{2\pi }{a}\Big)\int_{
    1/n}^{\infty}\frac{dt}{t^2}
    =2a^s\omega \Big(\frac{2\pi }{a}\Big)n$$
    $$\leqslant 2a^s\Big(\frac{2\pi n}{a}+1\Big)\omega \Big(\frac{1}{n}
    \Big)n<
    8a^{s-1}\pi n^2\omega \Big(\frac{1}{n}\Big)\leqslant
    8\pi n^{s+1}\omega \Big(\frac{1}{n}\Big),\ i=0,1.
\end{equation}
The estimate
\begin{equation}\label{7/08-7}
    \int_{-\infty}^{-1/n}\delta (t)\frac{1}{t}
    \int_{0}^{a}u^s
    \sin\Big(ut+\frac{i\pi }{2}\Big)\,du\,dt=O(1)n^{s+1}
    \omega (1/n),\ \ \ i=0,1
\end{equation}
is similarly proved. Comparing relations (\ref{7/08-6}), (\ref{7/08-7}) and
(\ref{11/08-10}), we obtain (\ref{7/08-3}). Lemma \hyperlink{18.08.10-16:25:32}{2} is
proved.\hspace{\stretch{1}}$\blacksquare$

Applying Lemma \hyperlink{18.08.10-16:25:32}{2} to the integral in (\ref{23.12.09-17:34:42})
and, at the same time, to $J_{3,n},$ we have
\begin{equation}\label{24.12.09-09:53:22}
    J_{1,n}=O(1)\psi (n)\omega (1/n),
\end{equation}
\begin{equation}\label{24.12.09-10:01:06}
    J_{3,n}=O(1)\psi (n)\omega (1/n).
\end{equation}
In the monograph \cite[pp. 212, 216, see relations (4.26$'$) and (4.42),
(4.45), (4.46)]{Stepanets_2002} it is shown, that
\begin{equation}\label{24.12.09-09:57:24}
    J_{2,n}=O(1)\psi (n)\omega (1/n)\ \ \ \forall \psi \in \mathfrak{M}_0
\end{equation}
and
\begin{equation}\label{24.12.09-09:58:05}
    J_{4,n}=O(1)\psi (n)\omega (1/n)\ \ \ \forall \psi \in
    \mathfrak{M}_0',\ \ \beta \neq 2l,\ \ l\in \mathbb{Z}.
\end{equation}
Since $|\tau _n(u)|\leqslant \psi (n),$ $u\in [0,1],$ it is clear that
\begin{equation}\label{24.12.09-10:03:13}
    J_{5,n}=O(1)\psi (n)\omega (1/n).
\end{equation}
Comparing (\ref{23.12.09-17:15:29}),
(\ref{24.12.09-09:53:22})--(\ref{24.12.09-10:03:13}), we arrive at
(\ref{23.12.09-17:29:28}). Then from (\ref{23.12.09-17:11:58}) for any
function $\varphi \in H_\omega ^0$ and $n\in \mathbb{N}$, we obtain
\begin{equation}\label{24.12.09-10:08:51}
    \int_{-\infty}^{\infty}\delta\Big(\frac{t}{n},\varphi\Big)
    \widehat{\tau }_n(t)\,dt$$
    $$=-\frac{\sin\frac{\beta \pi }{2}}{\pi }\int_{|
    t|\leqslant 1}\delta\Big(\frac{t}{n},\varphi\Big)
    \int_{1}^{\infty}\psi (nu)\sin ut\,du\,dt+O(1)\psi (n)\omega (1/n)$$
    $$=-\frac{\sin\frac{\beta \pi }{2}}{\pi }\int_{
    0}^{1}\bigg(\delta\Big(\frac{t}{n},\varphi\Big) -
    \delta\Big(-\frac{t}{n},\varphi\Big) \bigg)
    \int_{1}^{\infty}\psi (nu)\sin ut\,du\,dt$$
    $$+O(1)\psi (n)\omega (1/n),\ \ \ \psi \in\mathfrak{M}_0',\ \ \beta
    \in \mathbb{R}.
\end{equation}
Since
    \begin{equation}\label{24.12.09-15:18:57}
    \int_{1}^{\infty}\psi (nu)\sin ut\,du>0,\ \ \ t\in(0,1],\ \ \ \psi \in
    \mathfrak{M}',\ \ \beta \neq 2l,\ l\in\mathbb{Z}
    \end{equation}
(see, e.g., \cite[p. 143]{STEANETS-RUKASOV-CHAICHENKO}) and
\begin{equation}\label{24.12.09-10:19:15}
    \int_{0}^{1}\omega \Big(\frac{2t}{n}\Big)\int_{
    1}^{\infty}\psi (nu)\sin ut\,du\,dt$$
    $$=\int_{0}^{1/n}\psi \Big(
    \frac{1}{t}\Big)\frac{\omega (t)}{t}\,dt+O(1)\psi (n)\omega (1/n),\ \
    \ \psi \in \mathfrak{M}_0',\ \ \beta \neq 2l,\ l\in\mathbb{Z},
\end{equation}
(see \cite[p. 632]{SERDUK-OVSIY UMJ}), from (\ref{23.12.09-15:42:55}) and
(\ref{24.12.09-10:08:51}) we obtain (\ref{23.12.09-13:55:31}). Putting
together inequalities (\ref{23.12.09-13:55:31}) and
(\ref{23.12.09-13:55:02}) we find a required estimate for quantity
(\ref{13.12.09-00:56:41}):
\begin{equation}\label{24.12.09-10:36:36}
    E_n(C^{\psi}_{\beta}H_\omega)\leqslant
    \frac{1}{\pi }\big|\sin\frac{\beta \pi }{2}\big|\int_{0}^{1/n}\psi \Big(
    \frac{1}{t}\Big)\frac{\omega (t)}{t}\,dt$$
    $$+O(1)\psi (n)\omega (1/n),\ \
    \psi \in \mathfrak{M}_0',\ \ \beta \in\mathbb{R}.
\end{equation}

\textbf{2.} Let us find a lower bound for $E_n(C^{\psi}_{\beta}H_\omega).$

Let
 $\varphi _{n}(t)$ be an odd $2\pi/n$-periodic function defined on
 $[0,\pi /n]$ by the equalities
    $$\varphi _{n}(t)=
  \begin{cases}
    \frac{c_\omega }{2}\omega (2t), & t\in[0,\pi /2n], \\
    \frac{c_\omega }{2}\omega (\frac{2\pi}{n}-2t ), & t\in[\pi /2n,\pi /n],
  \end{cases}
    $$
where $c_\omega =1$ if $\omega (t)$ is a convex upwards modulus of
continuity and $c_\omega =2/3$ otherwise. As shown in \cite[pp.
83--85]{Stepanets_1987}
 if $\omega (t)$ is an arbitrary modulus
of continuity, then
    $$|\varphi _{n}(t')-\varphi _{n}(t'')
    |\leqslant \omega (|t'-t''|),\ \ t',t''\in[-\pi /2n,\pi /2n].$$
This implies that
    $$|\varphi _{n}(t')-\varphi _{n}(t'')|\leqslant
    \omega (|t'-t''|),\ \ t',t''\in \mathbb{R},$$
and, hence, $\varphi _{n}\in H_\omega .$ We denote by $f^*(\cdot)$ the
function from the set $C^{\psi}_{\beta}H_\omega$, \mbox{$\psi \in
\mathfrak{M}',$} whose $(\psi ,\beta )$-derivative $f^{*\psi
}_{\,\,\,\beta} (t)$ coincides with the function $\varphi _{n}(t)$ on a
period. By relations (\ref{24.12.09-14:19:33}), such a function
$f^*(\cdot)$ exists.

In virtue of formula (3.4) from the book \cite[Chap. 2, Subsec.
3.1]{Stepanets_1987} the following equality holds for any $f\in
C^{\psi}_{\beta}H_\omega$, $\psi \in \mathfrak{M}'$:
\begin{equation}\label{11.12.09-14:44:01}
    f(x)-U_{n-1}(f;x;\Lambda)$$
    $$=\frac{1}{\pi }\int_{-\pi }^{\pi }f^\psi _\beta (x+t)\bigg(
    \sum\limits_{k=1}^{\infty}\psi (k)\cos\Big(kt+\frac{\beta \pi }{2}
    \Big)$$
    $$-\sum\limits_{k=1}^{n-1}\lambda _k^{(n)}\psi (k)\cos\Big(kt+\frac{\beta \pi }{2}
    \Big)\bigg)\,dt,\ \ x\in \mathbb{R},\ \ n\in \mathbb{N},
\end{equation}
where $U_{n-1}(f;x;\Lambda)$ is a trigonometric polynomial of the form
(\ref{24.12.09-14:09:14}), such that \mbox{$\lambda _n^{(n)}=0.$} Since
function $\varphi _{n}(t)$ is odd $2\pi /n$-periodic, the equalities
\begin{equation}\label{11.12.09-14:57:49}
    \int_{-\pi }^{\pi }\varphi _{n}(t)\sin kt\,dt=0,\ \ \
    k=1,2,\ldots,n-1,\ \ n\geqslant 2
\end{equation}
(see, e.g., \cite[p. 159]{KORNEICHUK-5}) and
    $$\varphi _{n}\Big(\frac{i\pi }{n}+t\Big)=
    (-1)^i\varphi _{n}(t),\ \ \ i\in\mathbb{Z},$$
hold. Then, using relation (\ref{11.12.09-14:44:01}) for $f^*(\cdot)$, we
obtain
\begin{equation}\label{11.12.09-15:05:33}
    f^*\big(\frac{i\pi }{n}\big)-U_{n-1}
    \big(f^*;\frac{i\pi }{n};\Lambda\big)$$
    $$=\frac{(-1)^i}{\pi }\int_{-\pi }^{\pi }\varphi _{n}(t)\bigg(\sum\limits_{k=1}^{\infty}\psi (k)\cos\Big(kt+\frac{\beta \pi }{2}
    \Big)-\sum\limits_{k=1}^{n-1}\lambda _k^{(n)}\psi (k)\cos\Big(kt+\frac{\beta \pi }{2}
    \Big)\bigg)\,dt$$
    $$=\frac{(-1)^i}{\pi }\int_{-\pi }^{\pi }\varphi _{n}(t)\sum\limits_{k=1}^{\infty}\psi (k)\cos\Big(kt+\frac{\beta \pi }{2}
    \Big)\,dt$$
    $$=\frac{(-1)^i}{\pi }\sin\frac{\beta \pi }{2}\int_{-\pi }^{
    \pi }\varphi _{n}(t)\sum\limits_{k=n}^{\infty}
    \psi (k)\sin kt\,dt,\ \ i\in\mathbb{Z},\ \ n=2,3,\ldots
\end{equation}

It is obvious from this that there exist $2n$ points
\mbox{$t_i=\frac{i\pi}{n}$,} \mbox{$i=0,1,\ldots,2n-1,$} on the period
$[0,2\pi )$ at which the difference
    $$f^*(x)-U_{n-1}
    (f^*;x;\Lambda)$$
takes values with alternating signs. Then by the de la Vall\'{e}e Poussin
theorem \cite{Valle'e Poussin} (see also \cite[p. 312]{Stepanets_1987},
\cite[p. 491]{Stepanets_2002}), we find
    \begin{equation}\label{11.12.09-15:32:52}
        E_n(f^*)\geqslant \frac{1}{\pi }\Big|\sin\frac{\beta \pi }{2}
        \int_{-\pi }^{
    \pi }\varphi _{n}(t)\sum\limits_{k=n}^{\infty}
    \psi (k)\sin kt\,dt\Big|,\ \ \psi \in \mathfrak{M}',
    \end{equation}
where
    $$E_n(f^*)=\inf\limits_{t_{n-1}}{\|
    f^*(\cdot)-t_{n-1}(\cdot)\|_C},\ \ \ n\in \mathbb{N}.$$
From (\ref{11.12.09-15:05:33}) and (\ref{11.12.09-15:32:52}) it follows,
in particular, that
\begin{equation}\label{11.12.09-15:52:00}
    E_n(f^*)\geqslant |f^*(0)-U_{n-1}
    (f^*;0;\Lambda)|,\ \ n=2,3,\ldots
\end{equation}
Inequality (\ref{11.12.09-15:52:00}) is satisfied for triangular matrix
$\Lambda =\|\lambda _k^{(n)}\|,$ $k=\overline{1,n}$, such that $\lambda
_n^{(n)}=0$. Let's define its remaining elements in the following way:
    $$\lambda _k^{(n)}=\lambda ^\psi (k/n),\ \ k=\overline{1,n-1},\ \ n\in \mathbb{N},$$
where $\lambda ^\psi (\cdot)$ is defined by (\ref{24.12.09-14:52:56}).
Since in this case
    $$U_{n-1}(f^*;0;\Lambda )=U_{n-1}^\psi (f^*;0),$$
then from (\ref{11.12.09-15:52:00}) we obtain, taking the inequality
$E_n(C^{\psi}_{\beta}H_\omega)\geqslant E_n(f^*)$ into account,
\begin{equation}\label{24.12.09-14:58:23}
    E_n(C^{\psi}_{\beta}H_\omega)\geqslant |f^*(0)-U_{n-1}^\psi
    (f^*;0)|,\ \ n=2,3,\ldots,\  \psi \in \mathfrak{M}'.
\end{equation}
By virtue of (\ref{23.12.09-15:37:11}) and (\ref{24.12.09-10:08:51})
\begin{equation}\label{24.12.09-15:02:39}
    f^*(0)-U_{n-1}^\psi (f^*;0)=\int_{-\infty}^{\infty}
    \bigg(f^{*\psi
    }_{\,\,\,\beta} \Big(\frac{t}{n}\Big)-f^{*\psi }_{\,\,\,\beta}
    (0)\bigg)\widehat{\tau }_n(t)\,dt$$
    $$=-\frac{\sin\frac{\beta \pi }{2}}{\pi }\int_{0}^{1}\bigg(
    \varphi _n\Big(\frac{t}{n}\Big)-\varphi _n
    \Big(-\frac{t}{n}\Big)\bigg) \int_{1}^{\infty} \psi (nu)\sin
    ut\,du\,dt$$
    $$+
    O(1)\psi (n)\omega (1/n)$$
    $$=-c_\omega \frac{\sin\frac{\beta \pi }{2}}{\pi }\int_{0}^{1}
    \omega \Big(\frac{2t}{n}\Big)\int_{
    1}^{\infty}\psi (nu)\sin ut\,du\,dt+O(1)\psi (n)\omega (1/n),\ \psi
    \in \mathfrak{M}_0'.
\end{equation}
Combining (\ref{24.12.09-15:18:57}), (\ref{24.12.09-10:19:15}),
(\ref{24.12.09-14:58:23}) and (\ref{24.12.09-15:02:39}), we arrive at the
desired estimate
\begin{equation}\label{24.12.09-15:19:55}
    E_n(C^{\psi}_{\beta}H_\omega)\geqslant
    \frac{c_\omega }{\pi }\big|\sin\frac{\beta \pi }{2}\big|\int_{0}^{1/n}\psi \Big(
    \frac{1}{t}\Big)\frac{\omega (t)}{t}\,dt$$$$
    +O(1)\psi (n)\omega (1/n),\ \
    \psi \in\mathfrak{M}_0',\ \ \beta \neq 2l,\ l\in \mathbb{Z}.
\end{equation}
From (\ref{24.12.09-10:36:36}) and (\ref{24.12.09-15:19:55}) we obtain formula
(\ref{19.12.09-20:55:52}). Theorem \hyperlink{14.04.11-19:15:10}{1} is
proved.\hspace{\stretch{1}}$\blacksquare$

\emph{Contact information}: \href{http://www.imath.kiev.ua/~funct}{Department of the Theory
of Functions}, Institute of Mathematics of Ukrainian National Academy of Sciences, 3,
Tereshenkivska st., 01601, Kyiv, Ukraine \vskip 0.2 cm \emph{E-mail}:
\href{mailto:serdyuk@imath.kiev.ua}{serdyuk@imath.kiev.ua},
\href{mailto:ievgen.ovsii@gmail.com}{ievgen.ovsii@gmail.com}

\end{document}